\documentclass[11pt]{amsart}  
\usepackage{amssymb, amsthm, graphicx, epsfig, graphs, verbatim, amscd,
psfrag, doublespace, enumerate }  \textwidth=15cm 
\newtheorem{thm}{Theorem}[section] 
 
\newtheorem{lem}[thm]{Lemma} 
\newtheorem{prop}[thm]{Proposition} 
 
\theoremstyle{definition} 
\newtheorem{defn}[thm]{Definition}

\newtheorem{rem}[thm]{Remark}

\numberwithin{equation}{section}

\newcommand{\ep}{\varepsilon}

\newcommand{\si}{\sigma}

\newcommand{\x}{\times}

\newcommand{\Z}{\mathbb Z} 
\newcommand{\N}{\mathbb N} 
 
\newcommand{\R}{\mathbb R}

\newcommand{\del}{\partial}

\newcommand{\co}{\colon\thinspace}

\begin{document} 
\mathsurround=1pt  
\title[Morse functions on unoriented surfaces]{Cobordism group of Morse functions on unoriented surfaces}


\author{Boldizs\'{a}r Kalm\'{a}r} 
\address{E\"{o}tv\"{o}s Lor\'{a}nd University, P\'{a}zm\'{a}ny P\'{e}ter s\'{e}t\'{a}ny 1/c.  
H-1117 Budapest, Hungary} 
 
\email{kalmbold@cs.elte.hu}

\begin{abstract}

Ikegami and Saeki have proved that the cobordism group of Morse functions on 
oriented surfaces is an infinite cyclic group.  Their method is applicable 
with a little modification to the computation of the cobordism group of Morse
functions on unoriented  surfaces. We prove that this group is isomorphic to
the direct sum of the infinite cyclic group and the finite group of order two. 
%$\Z \oplus \Z_2$.  
 
\end{abstract}

\maketitle 

\begin{spacing}{1.3}
 
\section{Introduction}\label{s:intro}

Ikegami and Saeki \cite{IS} initiated the study of the oriented cobordism 
group of Morse functions on oriented  
surfaces. 
Summarizing in an informal language, two Morse functions $f_0 \co M_0 \to \R$ and $f_1 \co M_1 \to \R$ 
on closed surfaces $M_0$ and $M_1$ are  
{\it cobordant}, if there exists a smooth map 
$F \co X \to \R \times [0,1]$ with only fold singularities (i.e., a fold map;
 fold singularities are the higher dimensional
analogues of the non-degenerate critical points of Morse functions; for a precise definition
of a fold map, see \S 2)
from a compact 3-manifold $X$ with 
$\del X = M_0 \amalg M_1$, and $F$ restricted to $\del X$ 
is identified with $f_0 \amalg f_1$, where $\amalg$ denotes the disjoint union 
(while the normal direction is mapped onto the normal
direction isomorphically). 
Such an $F$ is a {\it cobordism} between $f_0$ and $f_1$.
One can define the {\it cobordism group} of Morse functions on surfaces, which describes the structure of 
Morse functions on surfaces from the viewpoint of the cobordism  relation.
% For cobordisms in general 
%see for example \cite{Sw} Chapter 12.
 Ikegami and Saeki \cite{IS}  have transformed the problem of Morse functions and their cobordisms into 
 a problem of maps from graphs to $\R$ and maps from 2-dimensional 
polyhedra to $\R \x [0,1]$. 
This transformation has been realized by using the Stein factorization: 
the Stein factorization of a cobordism between two Morse functions
provides a cobordism  in an abstract sense between the Stein factorizations of
the Morse functions (called also the Reeb functions). 

The computation of the oriented cobordism group of Morse functions on oriented surfaces 
\cite{IS} involves two steps:  
\begin{itemize} 
\item[(1)] 
computing the abstract cobordism 
group of abstract Reeb functions, which are maps of graphs to $\R$ like Reeb
functions, where a cobordism between two such maps is a map from a compact
2-dimensional polyhedron to $\R  \x [0,1]$ like the Stein factorization of a
cobordism between two Morse functions,
\item[(2)]  
showing that the abstract cobordism group is isomorphic to the cobordism group 
of Morse functions.  
\end{itemize} 
 This method enables us to solve the problem in 
a purely combinatorial context.

%The relationships are illustrated in the next
%diagram:\\   
% 
%\[ 
%\begin{graph}(6.3,4) 
%\opaquetextfalse 
%\graphlinecolour{1}\grapharrowtype{2} 
%\textnode {A}(1,4){Morse function $f_1$} 
%\textnode {B}(4.5,4){\footnotesize{cobordism by $F$}} 
%\textnode {C}(8,4){Morse function $f_2$} 
%\textnode {D}(1,2.5){Reeb function ${\bar f}_1$} 
%\textnode {E}(4.5, 2.5){\footnotesize{cobordism by $\bar F$}} 
%\textnode {F}(8, 2.5){Reeb function ${\bar f}_2$} 
%\textnode {G}(1, 1){an abstract Reeb function} 
%\textnode {H}(4.5, 1){\footnotesize{cobordism}} 
%\textnode {I}(8, 1){an abstract Reeb function} 
%\diredge {A}{D}[\graphlinecolour{0}] 
%\diredge {B}{E}[\graphlinecolour{0}] 
%\diredge {C}{F}[\graphlinecolour{0}] 
%\edge {D}{G}[\graphlinecolour{0}] 
%\edge {E}{H}[\graphlinecolour{0}] 
%\edge {F}{I}[\graphlinecolour{0}] 
%\freetext (-0.2,3.3){\tiny{Stein factorization}} 
%\freetext (3.3,3.3){\tiny{Stein factorization}} 
%\freetext (6.8,3.3){\tiny{Stein factorization}} 
%\freetext (-2,4.4){\tiny{cobordism group of}} 
%\freetext (-2,4.1){\tiny{Morse functions}} 
%\freetext (-2.2,1){\tiny{abstract}} 
%\freetext (-2.2,0.7){\tiny{cobordism}} 
%\freetext (-2.2,0.4){\tiny{group}} 
%\end{graph} 
%\] 
% 
 		
More precisely step (1) consists of the following. 
A cobordism between two Morse functions is a generic map of a compact
3-dimensional manifold into $\R \x [0,1]$. The local structure of the
Stein factorizations of such maps is completely described in \cite{Kush, Lev}.
Using this description, Ikegami and Saeki \cite{IS} gave an alternative
description of an abstract cobordism of Reeb  functions, which is more
algorithmic: if we change a function on a graph 
by a finite iteration of simple local moves together with a homotopy,
 we get a function cobordant to the  original function in the
abstract sense. Using this observation, Ikegami and Saeki \cite{IS} gave 
representatives of the elements of the abstract cobordism group, and proved 
that it is isomorphic to $\Z$. 
In step (2) \cite{IS} realizes  an abstract  
cobordism between Reeb functions of two Morse functions  as the Stein
factorization of a fold map of an oriented 3-manifold.  

In this paper, we follow
the same line of arguments in order to compute the cobordism group of Morse
functions\footnote{The author was informed, when this paper had already been
written, that  Ikegami \cite{Ik} had computed the analogous group in
each dimension. His method is quite different and more sophisticated than our
approach.} on unoriented surfaces. However, we have to modify the method of
\cite{IS}: in step (1), when  we want to compute the abstract
cobordism group, we have to consider functions and maps coming from the Stein
factorizations of Morse functions on {\it nonorientable} surfaces and those of
fold maps of {\it nonorientable} 3-manifolds, respectively. Levine \cite{Lev}
describes these maps completely. With  these new cases, we get an equivalent
algorithmic description of an abstract cobordism,  which differs from \cite{IS}
only in some new moves on a few new types of  Reeb graphs.  It turns out that
the method of the computation of the  abstract cobordism group given in
\cite{IS} is applicable in our case as well: we give a set  of representatives
of the elements of this group by using the algorithm  of \cite{IS}, and 
 we show
that the group is isomorphic to $\Z \oplus \Z_2$.    Step (2) also  needs some
modifications. On the one hand this step becomes simpler since we are 
allowed to use nonorientable 3-manifolds for the realization of abstract
cobordisms   between Reeb functions. On the other hand this step becomes more
complicated, since the   method used in \cite{IS} cannot be applied here.
An abstract cobordism between two given Reeb functions is a map from a compact
2-dimensional polyhedron to $\R  \x [0,1]$. 
In order to realize the neighbourhood of an appropriate 1-skeleton  of this polyhedron as a Stein 
factorization of a fold map of a 3-manifold $X$ with boundary, one can apply the 
method described 
in \cite{IS}, but this realization is not a cobordism between the two given Reeb 
functions. To obtain a cobordism we attach  2-disk bundles and their fold maps to the 
 appropriate part of the boundary of $X$, as will be described in \S 4.

The author would like to thank Andr\'{a}s Stipsicz and Andr\'{a}s
Sz\H{u}cs for their help and encouragement.   

\section{Preliminaries}\label{s:prelim}

\begin{defn} 
 
Let $Q^q$ and $N^n$ be smooth manifolds of dimensions $q$ and $n$ 
respectively $(q \geq n)$. Let $p \in Q^q$ be a singular point of $f$.
A smooth map $f \co Q^{q} \to N^{n} $  has a {\it fold 
singularity} at the singular point $p$, if we can write $f$ in some local coordinates at $p$  
and $f(p)$ in the form 
\[  
f(x_1,\ldots,x_q)=(x_1,\ldots,x_{n-1},\pm x_n^2 \pm \cdots \pm x_q^2).
\] 
A smooth map $f \co Q^{q} \to N^{n}$ is called a {\it fold map}, if $f$ has only 
fold singularities.
A smooth function on a manifold which has only fold singularities is called a
{\it Morse function}.

\end{defn} 
Let $\amalg$ denote the disjoint union.
\begin{defn} 
 
Two Morse functions $f_0 \co M_0 \to \R$ and $f_1 \co M_1 \to \R$ 
on closed surfaces $M_0$ and $M_1$ are  
{\it cobordant}, if there exists a smooth map  
$F \co X \to \R \times [0,1]$ with only fold singularities from a compact 3-manifold $X$ with
$\del X = M_0 \amalg M_1$ such that  ${F \mid}_{M_0 \x [0,\ep)}=f_0 \x
$id$_{[0,\ep)}$ and ${F \mid}_{M_1 \x (1-\ep,1]}=f_1 \x $id$_{(1-\ep,1]}$, where $M_0 \x [0,\ep)$
 and $M_1 \x (1-\ep,1]$ are small collar neighbourhoods of $\del X$ with the
identifications $M_0 = M_0 \x \{0\}$, $M_1 = M_1 \x \{1\}$, and id$_A (a) = a$ $(a \in A)$
for an arbitrary set $A$.

We call  the map $F$ a {\it cobordism} between $f_0$ and $f_1$.

\end{defn} 
 This clearly defines an equivalence relation. 

\begin{defn} 
We can define a group operation on the set
of the cobordism classes of Morse functions as follows.
If $[f_0 \co M_0 \to \R]$ and 
$[f_1 \co M_1 \to \R]$ are cobordism classes of Morse functions, then the {\it sum} of $[f_0]$ and $[f_1]$ is
the cobordism class represented by the Morse function $f_0 \amalg f_1 \co M_0 \amalg M_1 \to \R$. 

It is easy to show that the cobordism class of the Morse function $f_0 \amalg f_1$ does not
depend on the choice of the representatives $f_0$ and $f_1$.
\end{defn}
This group operation is clearly commutative, so we have an abelian group, which we denote by
${\mathcal Cob}_f(2, -1)$, where the lower index ``$f$'' comes from ``fold map''.

\begin{rem}
One can define the analogous cobordism groups ${\mathcal Cob}_f(q, k)$ ($q \in \N$, $k \in \Z$, $-q \leq k \leq 0$)
of fold maps $f \co Q^q \to \R^{q+k}$ from closed $q$-dimensional manifolds $Q^q$ to $\R^{q+k}$ as well.
The integer $k$ is called the {\it codimension} of the map $f$. 

For $k \geq 0$ there are many results concerning the cobordism groups ${\mathcal Cob}_{\tau}(q, k)$,
where $\tau$ is a set of singularity types, the elements of the group ${\mathcal Cob}_{\tau}(q, k)$ 
are cobordism classes of smooth maps with only singularities in $\tau$, and a cobordism between two 
such maps has only singularities in $\tau$. See for example \cite{Ko, RSz, Szucs1, Szucs2}. 

\end{rem}
 
The aim of this paper is to prove 
 
\begin{thm} \label{mainthm}
 
The cobordism group ${\mathcal Cob}_f(2, -1)$ is isomorphic to $\Z \oplus \Z_2 
$. 
 
\end{thm}

We use the notion of the Stein factorization of a smooth map 
$f\co Q^q \to N^n$, where $Q^q$ and $N^n$ are smooth manifolds 
of dimensions $q$ and $n$ respectively $(q \geq n)$. 
Two points  $p_1,p_2 \in Q^q$ are {\it equivalent} if and only if 
$p_1$ and $p_2$ lie on the same component of an $f$-fiber. 
Let $W_f$ denote the quotient space of $Q^q$ with respect
to this equivalence relation and $q_f \co W_f \to N^n$ the quotient map.
Then there exists a unique continuous map $\Bar{f} \co W_f \to N^n$ such that
$f = \Bar{f} \circ q_f$. The space $W_f$ or the factorization of the 
map $f$ into the composition of $q_f$ and $\Bar{f}$ is called the {\it Stein
factorization} of the map $f$. We call the map $\Bar{f}$  the  
{\it Stein factorization} of the map $f$ as well.

 In the case of a Morse function $f \co Q^q \to \R$ on a closed manifold $Q^q$, 
$W_f$ is a graph (since every regular fiber is a 1-codimensional submanifold
in $Q^q$ and the set of critical points of $f$ is finite and discrete), 
which is called the {\it Reeb graph}, and the 
function  
$\bar{f} \co W_f \to \R$ is called the {\it Reeb
function} of the Morse function $f$.

% (for example, see \cite{Fom}).    

We may assume that 
the Morse function $f$ is in general position, that is, every fiber contains at most 
one critical point, since a small perturbation of the function $f$ is a Morse
function cobordant 
to $f$.

\section{Abstract cobordism}\label{s:chap} 
 
\begin{figure}[ht]  
\begin{center}  
\psfrag{a}{(a)}
\psfrag{b}{(b)}
\psfrag{c}{(c)}
\psfrag{p}{$+1$}
\psfrag{m}{$-1$}
\epsfig{file=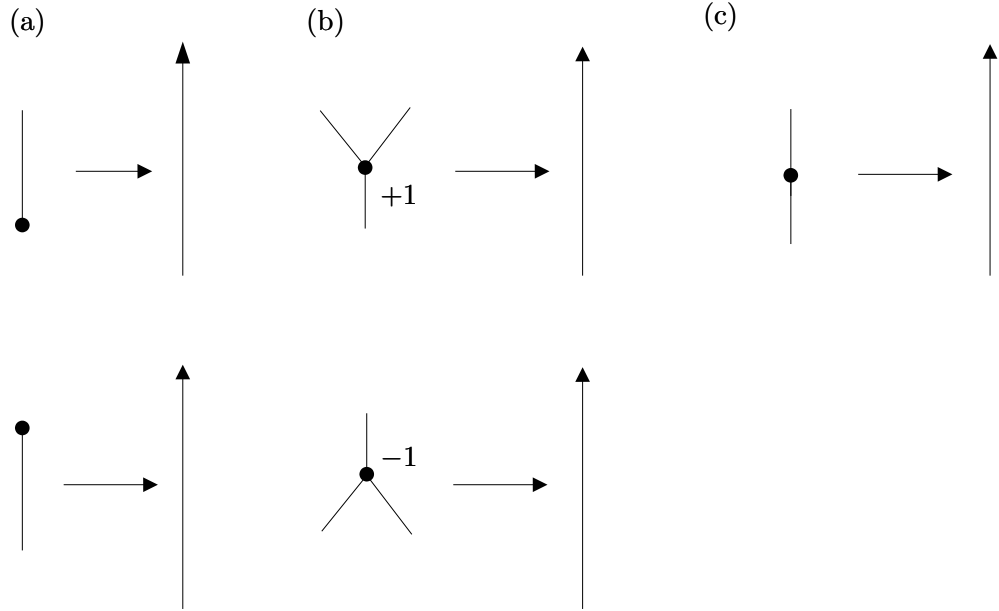, height=10cm}  
\end{center}  
\caption{The Reeb function of a Morse function  $f$ in a neighbourhood of
the $q_f$-image of a critical point}  
\label{f:functions}  
\end{figure}

Let $f \co M \to \R$ be a Morse function on a closed surface $M$.
One can see (for example, see \cite{Lev, Sa}) that the Reeb function of $f$ 
in a neighbourhood of  the $q_f$-image of a critical point
is equivalent to one of the functions as depicted
in Fig.~\ref{f:functions}.   
In case (a) the corresponding critical point of $f$ is a local minimum or 
a local maximum, in case (b) the corresponding critical point of $f$ has 
index 1 with sign $+1$ or $-1$, and a sign is associated to each 
vertex of degree three as in Fig.~\ref{f:functions}.
 Case (c) can
occur only on a nonorientable surface. The corresponding critical  
point of $f$ is also of index 1.  Note that the Reeb function  
restricted to an edge of the Reeb graph is always an embedding.

\begin{defn}

We say that a smooth map $F \co X \to \R^2$ from a closed
3-manifold $X$  has a {\it definite fold
singularity} at a singular point $p \in X$ if we can write $F$ in some local coordinates at $p$  
and $F(p)$ in the form 
\[  
F(x_1,x_2,x_3)=(x_1,x_2^2 + x_3^2),
\] 
and we say that $F$ has an {\it indefinite fold singularity} at a singular point $p \in X$ if
 we can write $F$ in some local coordinates at $p$  
and $F(p)$ in the form 
\[  
F(x_1,x_2,x_3)=(x_1,x_2^2 - x_3^2).
\]
%Let $S_F$ denote the set of singular points of $F$ in $X$. 
We say that $F$ has a {\it simple indefinite fold singularity}
at a singular point $p \in X$ if the map $F$ has an indefinite fold singularity at $p$ and 
each component of $F^{-1}(F(p))$ contains at most one singular point.
%$F^{-1}(F(p)) \cap S_F = \{p\}$. 
Furthermore, $F$ has a {\it nonsimple indefinite fold singularity} 
at $p$ if $F$ has an indefinite fold singularity at $p$ and 
there is a component of $F^{-1}(F(p))$ which contains at least two singular points.
%$F^{-1}(F(p)) \cap S_F \neq \{p\}$.

\end{defn}

\begin{figure}[ht]  
\begin{center}  
\psfrag{a}{(a)}
\psfrag{b}{(b)}
\psfrag{c}{(c)}
\psfrag{d}{(d)}
\psfrag{e}{(e)}
\psfrag{f}{(f)}
\psfrag{g}{(g)}
\psfrag{h}{(h)}
\psfrag{i}{(i)}
\epsfig{file=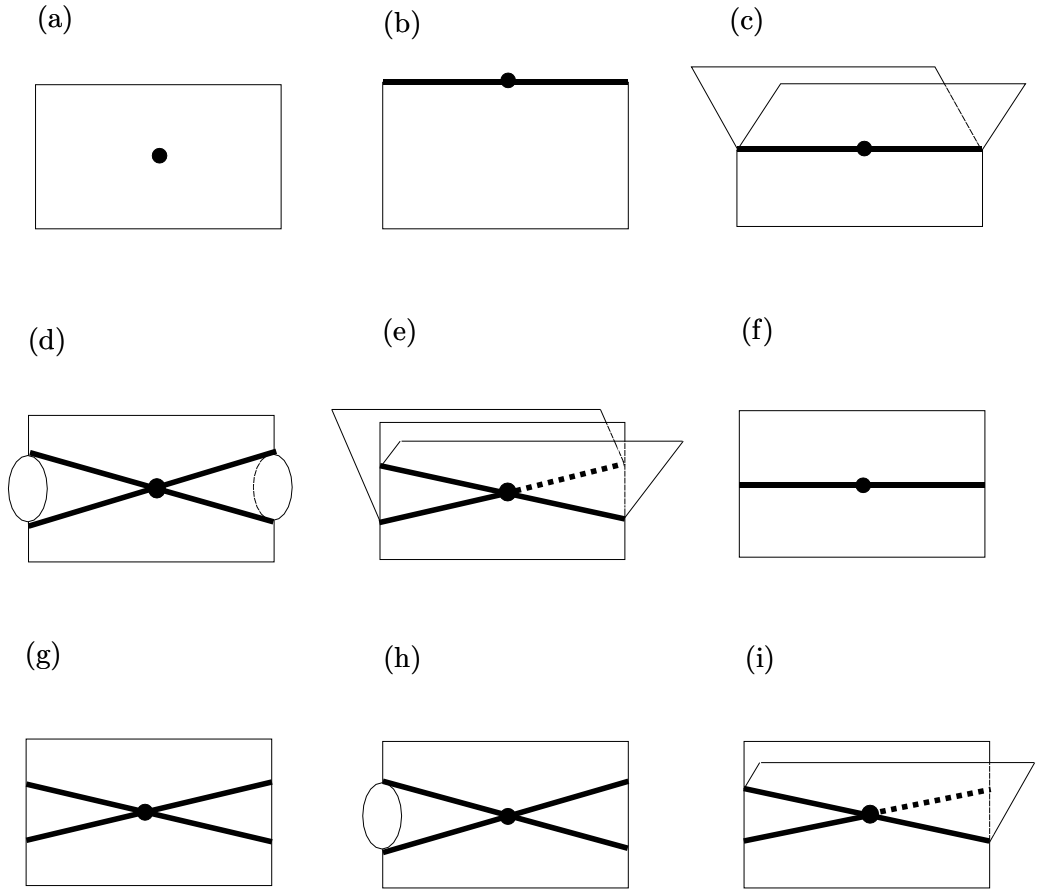, height=13cm}  
\end{center}  
\caption{The Stein factorization of a fold map $F$ in a 
neighbourhood of the $q_F$-image of a singular point}  
\label{f:maps} 
\end{figure}

Let $F \co X \to \R^2$ be a generic smooth map of a closed 3-manifold $X$
with only fold singularities.
By \cite{Lev} the Stein factorization $\bar{F} \co W_{F} \to \R^2$ of
 $F$  in a neighbourhood of 
the $q_F$-image of a singular point is equivalent to one of the maps as
depicted in Fig.~\ref{f:maps}, where the 
maps in question are projections to the plane of the paper. 
 The pictures (f)--(i) can occur only if $X$ is nonorientable. The lines in 
bold indicate the $q_F$-image of the set of singular points of $F$, which is a
1-dimensional submanifold of $X$ and is immersed into $\R^2$ by $F$. The dot in
bold indicates  the $q_F$-image of the singular point under study.  Case (a)
corresponds to a regular value, case (b) corresponds to the image of a
definite fold singularity, case (c) corresponds to the image of a simple 
indefinite fold singularity, cases (d) and (e) correspond to a nonsimple
indefinite fold singularity, and these exhaust all the possibilities 
when $X$ is orientable. Case (f) corresponds to the image of a simple
indefinite fold  singularity, and cases (g), (h), and (i) correspond to a
nonsimple indefinite fold singularity. Details can be found in \cite{IS} ,
\cite{Kush} and \cite{Lev}.        

Based on these observations, we give a definition of abstract Reeb functions
and their cobordism group analogous to \cite{IS} as follows.    
\begin{defn} \label{def1}
Let $G$ be a finite graph such that each of its vertices is of degree 1, 2 or 3. 
A continuous function $r \co G \to \R$ is said to be an {\it abstract Reeb function}
if $r$ is an embedding on each edge, and $r$ is equivalent in a neighbourhood 
of every vertex to one of the functions as depicted in Fig.~\ref{f:functions}.
\end{defn}   
\begin{defn} 
Let $r_i \co G_i \to \R$, $i = 0,1$, be abstract Reeb functions in the sense of 
Definition~\ref{def1}. We say that $r_0$ and $r_1$ are {\it cobordant} if there
exists a continuous map $R \co P \to \R \x [0,1]$ of a compact 2-dimensional 
polyhedron $P$ such that 
\begin{itemize}
\item[(i)]
$G_i = R^{-1}( \R \x \{ i \} )$, $i = 0,1$, are subcomplexes of $P$ with regular 
neighbourhoods of the forms $G_0 \x [0,\ep)$ and $G_1 \x (1 - \ep, 1]$, where
$G_i$ corresponds to $G_i \x \{ i \}$, $i = 0,1$,
\item[(ii)]
${R \mid}_{G_0 \x [0,\ep)}=r_0 \x $id$_{[0,\ep)}$ and ${R \mid}_{G_1 \x 
(1-\ep,1]}=r_1 \x $id$_{(1-\ep,1]}$, 
\item[(iii)]
in a neighbourhood of each point of $P \setminus (G_0 \cup G_1)$, the polyhedron 
$P$ and the map $R$ are equivalent to one of the pictures as depicted in Fig.~\ref{f:maps}.
\end{itemize}
Furthermore, we call the map $R \co P \to \R \x [0,1]$ a {\it cobordism} between
$r_0$ and $r_1$.
\end{defn} 
This clearly defines an equivalence relation.

\begin{defn}
We can define a group operation on the set
of the cobordism classes of abstract Reeb functions as follows.
If $[r_0 \co G_0 \to \R]$ and 
$[r_1 \co G_1 \to \R]$ are cobordism classes of abstract Reeb functions, then the {\it sum} of $[r_0]$ and $[r_1]$ is
the cobordism class represented by the abstract Reeb function $r_0 \amalg r_1 \co G_0 \amalg G_1 \to \R$. 

It is easy to show that the cobordism class of the abstract Reeb function $r_0 \amalg r_1$ does not
depend on the choice of the representatives $r_0$ and $r_1$.
\end{defn}

Let us denote the cobordism group of abstract Reeb functions by $\mathcal A$.

\begin{prop} \label{propo}
 
The cobordism group ${\mathcal A}$ is isomorphic to $\Z \oplus \Z_2$.

\end{prop} 
 
\begin{proof} 
By the facts explained in \cite{IS}, we see that two abstract Reeb functions
$r_0 \co G_0  \to \R$ and $r_1 \co G_1 \to \R$ are cobordant if and only if we
can deform  $r_0$ to $r_1$ by a finite iteration of the eleven local moves as
depicted in Fig.~\ref{f:moves} up to a homotopy in the space of abstract
Reeb functions.    
\begin{figure}[ht]  
\begin{center}  
\psfrag{a}{(a)}
\psfrag{b}{(b)}
\psfrag{c}{(c)}
\psfrag{d}{(d)}
\psfrag{e}{(e)}
\psfrag{f}{(f)}
\psfrag{g}{(g)}
\psfrag{h}{(h)}
\psfrag{i}{(i)}
\psfrag{j}{(j)}
\psfrag{k}{(k)}
\epsfig{file=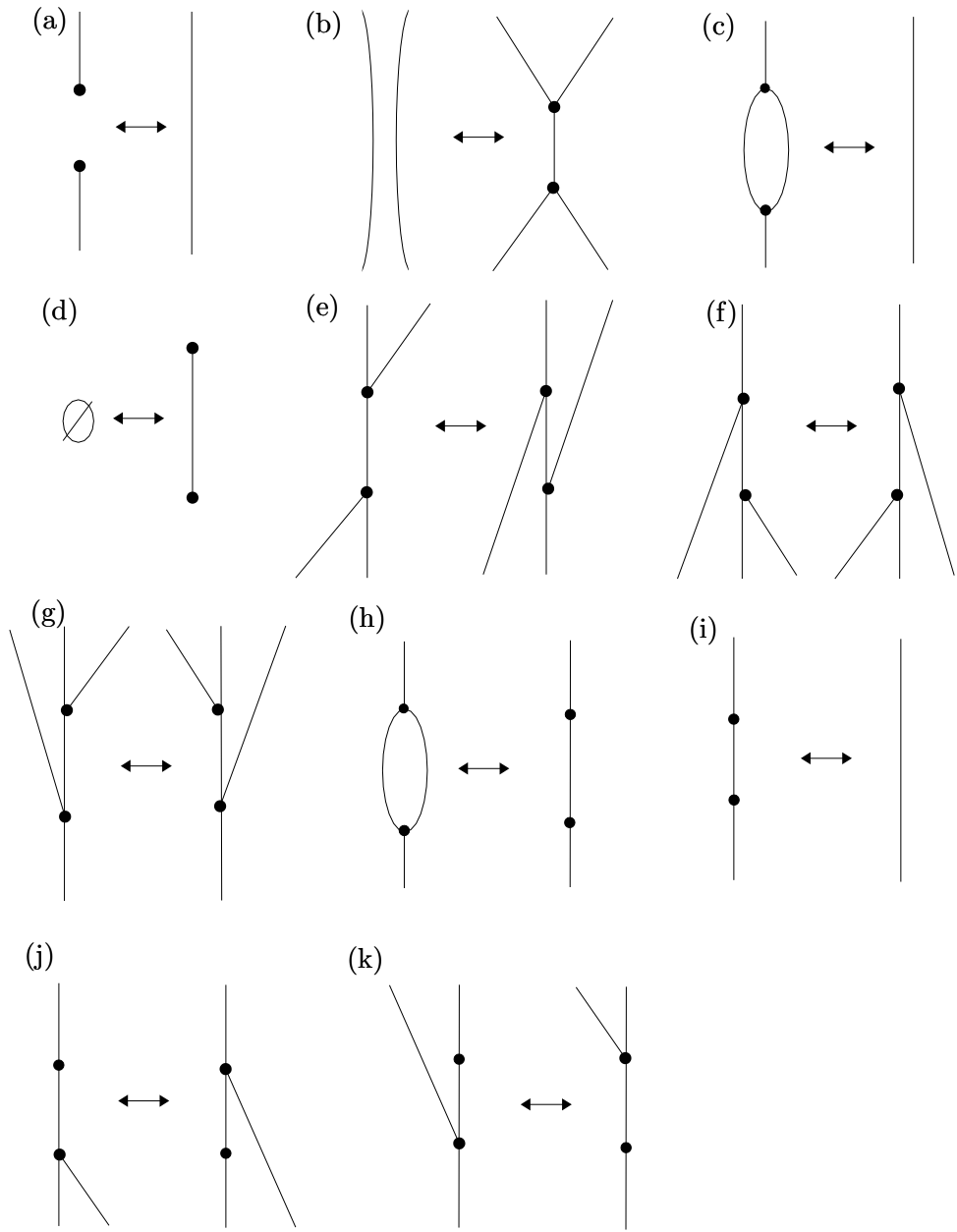, height=17cm}  
\end{center}  
\caption{Eleven local moves for abstract Reeb functions}  
\label{f:moves}  
\end{figure}   
 %In Figure~\ref{f:moves} every move corresponds to a case in  
%Figure~\ref{f:maps}, 
%when we pass through a singular point (details in \cite{IS}). More 
%precisely (b) and (c) in Figure~\ref{f:moves}  correspond to (c) in  
%Figure~\ref{f:maps}, 
%(f) and (g) to (e), (h) to (h), (i) to (f), (j) and (k) to (i). 
 
By \cite{IS} we can transform every abstract Reeb function whose domain
(source) graph has no  vertices of degree 2, using the first seven local moves
of Fig.~\ref{f:moves} together with a homotopy in the space of abstract Reeb
functions, into one of the abstract Reeb functions as depicted in
Fig.~\ref{f:reprori}, where  $n$ stands for the sum of the
signs over all vertices of degree 3. 
\begin{figure}[ht]  
\begin{center}
\psfrag{n}{$n$}
\psfrag{l}{$|n|$}
\psfrag{o}{$n=0$}
\psfrag{p}{$n > 0$}
\psfrag{m}{$n < 0$}
\epsfig{file=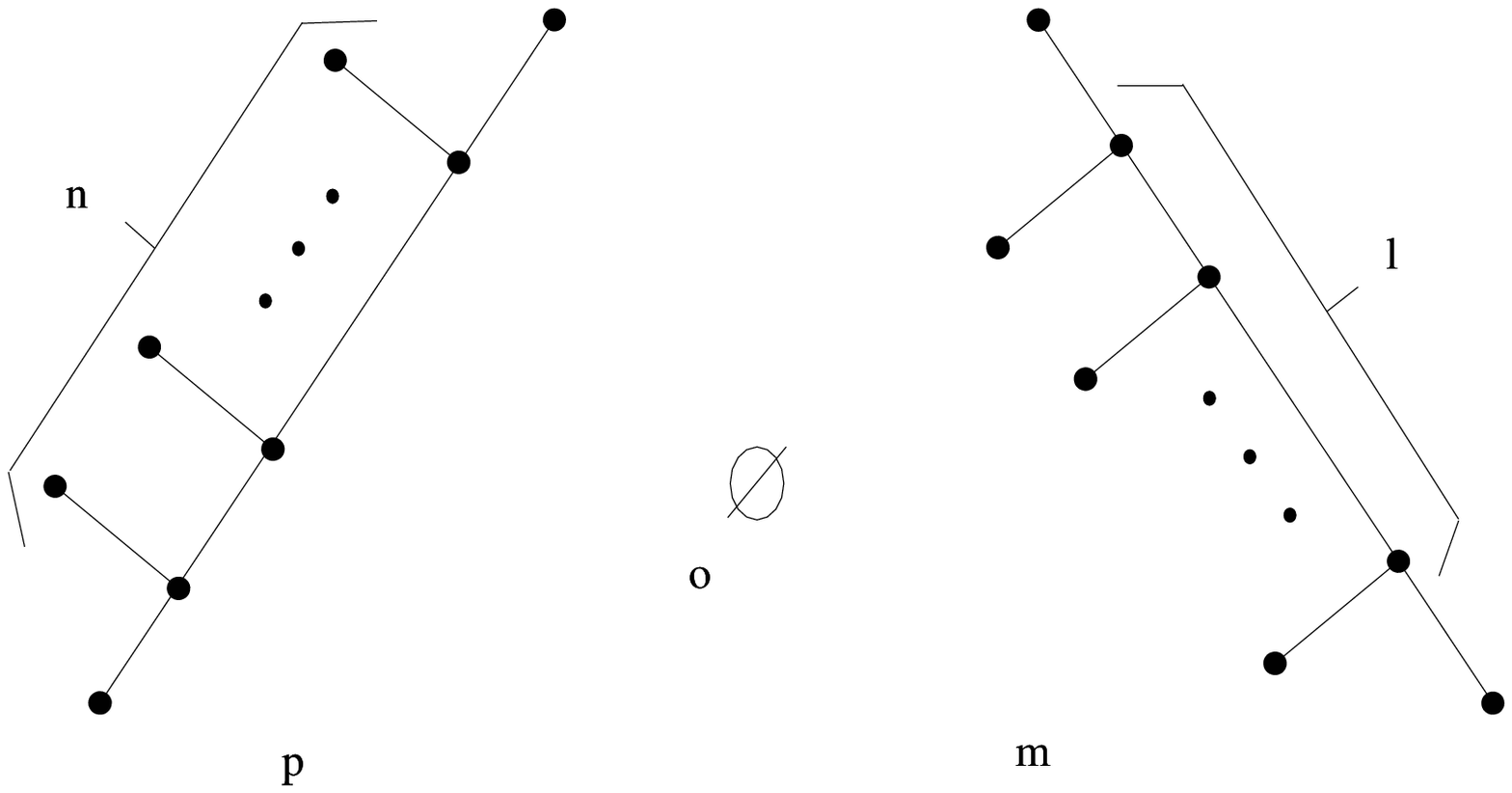, height=7cm}  
\end{center}  
\caption{First family of representatives of the elements of the cobordism
group ${\mathcal A}$}  
\label{f:reprori}   
\end{figure}   
 
\begin{figure}[ht]  
\begin{center} 
\psfrag{n}{$n$}
\psfrag{l}{$|n|$}
\psfrag{o}{$n=0$}
\psfrag{p}{$n > 0$}
\psfrag{m}{$n < 0$} 
\epsfig{file=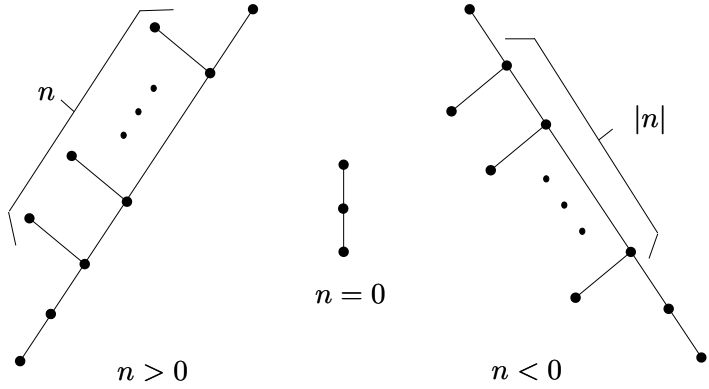, height=7cm}  
\end{center}  
\caption{Second family of representatives of the elements of the cobordism group ${\mathcal A}$} 
\label{f:reprnonori}   
\end{figure}

\begin{lem} \label{lemma}
Every abstract Reeb function is cobordant to one of the functions as depicted
in Figs.{\rm~\ref{f:reprori}} and {\rm \ref{f:reprnonori}}.
\end{lem} 
\begin{proof} 
With the help of the local moves (j) and (k) of Fig.~\ref{f:moves}, we can
apply the algorithm of \cite{IS} even if a graph has vertices of
degree 2, since we can lower any vertex of degree 2 through a vertex of
degree 3. Furthermore by the local move (i) of Fig.~\ref{f:moves}, we can
arrange so that the graph has at most one vertex of degree 2.
\end{proof} 

 Now let us define a map $\si \co {\mathcal A} \to \Z \oplus \Z_2$ from the
cobordism group of abstract Reeb functions as follows. 
 For a given abstract Reeb function $r$ put 
\[
\si(r) = (t(r), d(r)) \in \Z \oplus \Z_2,
\]
 where $t(r) \in \Z$ is the sum of the signs over all vertices of degree 3, and
$d(r) \in \Z_2$ is the number of vertices of degree 2 modulo two.
Then $\si(r)$ depends only on the cobordism class of an abstract Reeb
function $r$, since the local moves of Fig.~\ref{f:moves} or a homotopy in the
space of abstract Reeb functions do not change the values of $t$ nor $d$.
So $\si$ defines a well-defined map of the set of 
 cobordism 
classes, which is clearly a homomorphism, and by Lemma~\ref{lemma} it is an 
isomorphism. This completes the proof of Proposition~\ref{propo}.
\end{proof} 
 
\section{Relation between cobordism and abstract cobordism}\label{s:chap}

If we have a cobordism $F \co X \to \R \times [0,1]$ between two Morse 
functions $f_0 \co M_0 \to \R$ and $f_1 \co M_1 \to \R$ on closed surfaces,
then the Stein factorization 
$\bar{F} \co W_F \to \R \x [0,1]$ of the  possibly   perturbed $F$ provides a cobordism 
between the abstract Reeb functions $r_0 \co G_0 \to \R$ and $r_1 \co G_1 \to 
\R$ associated with the Morse functions $f_0$ and $f_1$ respectively, where
$r_j = {\Bar{f}}_j$ and $G_j = W_{f_j}$, $j = 0,1$ (for details, see
\cite{IS}). So the following definition makes sense.    

\begin{defn} 

 Let $\varphi \co {\mathcal Cob}_f(2,-1) \to \mathcal A$ be the map associating
to the cobordism  class of a Morse function $f$ the cobordism class of its
Reeb function $\Bar{f} \co W_f \to \R$.   This is clearly a well-defined 
homomorphism between the two cobordism groups 
 ${\mathcal Cob}_f(2,-1)$ and $\mathcal A$.

\end{defn} 
 
\begin{prop} \label{lenyeg}
 
The homomorphism $\varphi$ is an isomorphism. 
 
\end{prop} 
 
\begin{proof} 

The map $\varphi$ is clearly surjective. For the injectivity let $R \co P \to \R \x 
[0,1]$ be a cobordism between two abstract Reeb functions $r_0 \co G_0 \to 
\R$ and $r_1 \co G_1 \to \R$, and suppose that $r_0$ and $r_1$ correspond to 
Morse functions $f_0$ and $f_1$ on closed surfaces  $M_0$ and
$M_1$ respectively. 

It is sufficient to show that there exists a cobordism $F$ between the given Morse
functions $f_0$ and $f_1$. For this, it suffices to show that 
\begin{itemize}
\item[(1)]
there exists 
a cobordism $\tilde R \co 
\tilde P \to \R \x [0,1]$ between the two Reeb functions $r_0 \co G_0 \to 
\R$ and $r_1 \co G_1 \to \R$, and 
\item[(2)]
there exists a fold map
$\tilde F \co \tilde X \to \R \x [0,1]$ from a compact 3-manifold $\tilde X$ with boundary such
that $\tilde F$ is a composition of
a map $\tilde X \to \tilde P$ and the map $\tilde R$, and $\tilde R$ is identified with 
the Stein factorization 
of $\tilde F$. 
In this case we say that the map $\tilde F$ is {\it over} $\tilde R$.
\end{itemize}
Starting  in the same way as in \cite{IS}, we consider the
decomposition \[ P = N(Q) \cup N(V) \cup N(\Sigma) \cup T \]  as follows. We
put $Q = G_0 \cup G_1$, and let $N(Q)$ $(\cong (G_0 \x [0,\ep]) \amalg (G_1 \x [1-\ep,1]))$ 
be the regular neighbourhood of $Q$ in
$P$. Let $V$ be the set of the points of $P$ which have regular neighbourhoods
as in (d), (e), (g), (h), and  (i) of Fig.~\ref{f:maps}, and let $N(V)$ be 
its small regular
neighbourhood in $P$.  Let $\Sigma$ be the 
set of the points of $P$ which have regular neighbourhoods as in (b)--(i) of 
Fig.~\ref{f:maps}, and let $\tilde{N}(\Sigma)$ be its regular neighbourhood
in $P$.  Then $N( \Sigma )$ is the  closure of $\tilde{N} (
\Sigma ) \setminus ( N(Q) \cup N(V) )$.  Let $T$ be the closure of $P
\setminus ( N(Q) \cup N(V) \cup N(\Sigma))$. Note that $T$ is a compact
surface with boundary. 

By the construction of Mata-Lorenzo \cite{MaLo} we can
construct a 3-dimensional compact manifold $X$  with boundary and a map $F \co
X \to \R \x [0,1]$ over $R \mid_{N(Q) \cup N(V) \cup N(\Sigma)}$ such  that $\del X =
M_0 \amalg M_1 \amalg S$, ${F \mid}_{M_0 \x [0,\ep)}=f_0 \x  $id$_{[0,\ep)}$,
${F \mid}_{M_1 \x (1-\ep,1]}=f_1 \x $id$_{(1-\ep,1]}$, and $S$ is a
2-dimensional manifold fibered over $\del T$: 
in fact, we consider $(M_0 \x [0,\ep]) \amalg  (M_1 \x [1-\ep,1])$ over $N(Q)$,
and over $N(V) \cup N(\Sigma)$ we consider a manifold similar to that
constructed in \cite{IS}.

 It is not clear whether we can extend $X$ and $F$ like in  \cite{IS}
over $T$ so that the extension is a cobordism between the given Morse functions
$f_0$ and $f_1$, since now $S$ can be nonorientable.  {From}
\cite{Lev} $T$ is an  orientable surface with boundary and ${R \mid}_{T} \co T
\to \R \x [0,1]$ is an immersion. Over  $\del T$ we have an $S^1$-bundle
($=S$).     
Summarizing in an informal language, we have constructed a
3-manifold  $X$ with boundary $M_0 \amalg M_1 \amalg S$ and its fold map into
$\R \x [0,1]$ whose Stein  factorization gives a neighbourhood of the
1-skeleton of the 2-dimensional polyhedron $P$, where 1-skeleton means
$Q \cup \Sigma$.  It remained to extend this
construction over the 2-cells of $P$, that is over $T$. 
This was  easy in the oriented case
considered in \cite{IS}, since $S$ was a trivial  circle bundle over the
boundary of $T$ in $P$.  In the present case $S$ may contain a
nonorientable circle budle. Such a  nonorientable circle bundle is a fibration
of the Klein bottle over $S^1$.  This fibration can be extended to a fold map
of the (nonorientable) solid Klein bottle $S^1  {\x}_{{\Z}_2} D^2$ to the
annulus $S^1 \x [0,1]$ with only definite fold singularities, where the
boundary Klein bottle $S^1 {\x}_{{\Z}_2} S^1$ is mapped onto $S^1 \x \{1\}$ by
the above fibration, and the zero section $S^1 {\x}_{{\Z}_2} \{0\}$ coincides
with the set of the singular points and is mapped onto $S^1 \x
\{0\}$ diffeomorphically. (For this, it suffices to consider the function $(x,
y) \mapsto x^2 + y^2 $ on each 2-disk fiber.) Attaching such a solid Klein
bottle or its orientable version to $X$ along each component of $S$, we obtain
a compact 3-dimensional manifold and its fold map into $\R \x [0,1]$, which
gives a cobordism between the More functions $f_0$ and $f_1$. This completes
the proof of Proposition~\ref{lenyeg} and Theorem~\ref{mainthm}. \end{proof}

\begin{figure}[ht]  
\begin{center}  
\psfrag{a}{(a)}
\psfrag{b}{(b)}
\epsfig{file=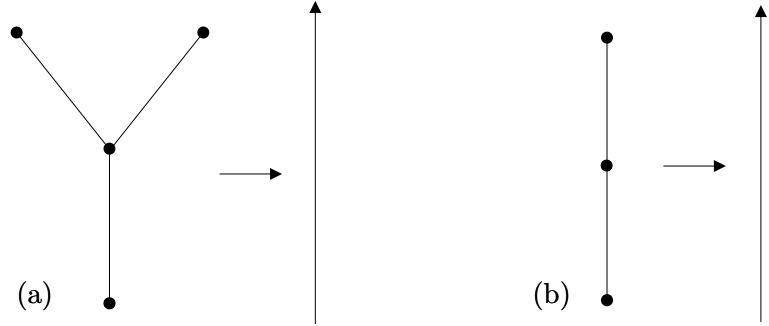, height=5cm}  
\end{center}  
\caption{Two generators of the abstract cobordism group}  
\label{f:generators}  
\end{figure}

\begin{rem} 
 
We see that the two functions as depicted in Fig.~\ref{f:generators} represent generators of the cobordism group of abstract Reeb 
functions. The function (a) has infinite order, while the function (b) has order 2. We can realize the function (a) 
as the Reeb function associated with a height function of a 2-sphere appropriately embedded in $\R^3$
 (see \cite{IS}), and the function (b) can be realized as the Reeb function associated with
a Morse function on $\R P^2$ with three critical points whose indices are 0, 1 and 2.

\end{rem}

\end{spacing} 
\end{document}